\documentclass{amsart}
\usepackage{palatino}
\usepackage{amsthm, amsmath}
\usepackage{amssymb} 
\usepackage{amscd}
\usepackage[hmargin=3.2cm, vmargin=2.8cm]{geometry}
\usepackage[breaklinks,bookmarksopen,bookmarksnumbered]{hyperref}
\usepackage[all]{xy}

\theoremstyle{plain}

\newtheorem*{thm*}{Theorem}

\newtheorem{prop}{Proposition}
\newtheorem{cor}{Corollary}
\theoremstyle{remark}
\newtheorem{rem}{Remark}

\newcommand\pr{\noindent\textit{Proof} : }

\newcommand\rond{\kern 1pt{\scriptstyle\circ}\kern 1pt}

\newcommand\Hom{\operatorname{Hom}}
\newcommand\im{\operatorname{Im}}
\newcommand\Ker{\operatorname{Ker}}
\newcommand\Coker{\operatorname{Coker}}

\newcommand\Tr{\operatorname{Tr}}

\newcommand\Z{\mathbb{Z}}
\newcommand\Q{\mathbb{Q}}
\newcommand\R{\mathbb{R}}

\renewcommand\P{\mathbb{P}}

\def\qfl#1{\buildrel {#1}\over {\longrightarrow}}

\newcommand\iso{\vbox{\hbox to .8cm{\hfill{$\scriptstyle\sim$}\hfill}
\nointerlineskip\hbox to .8cm{{\hfill$\longrightarrow $\hfill}} }}

\begin{document}
\title[On the second lower quotient of the fundamental group]{On the second lower quotient of the fundamental group}
\author[Arnaud Beauville]{Arnaud Beauville}
\address{Universit\'e de Nice Sophia Antipolis\\
Laboratoire J.-A. Dieudonn\'e\\
UMR 7351 du CNRS\\
Parc Valrose\\
F-06108 Nice cedex 2, France}
\email{arnaud.beauville@unice.fr}
 
 

\maketitle 
\section{Introduction}

Let  $ X $  be a connected topological space. The group $ G:= \pi_1(X)  $ admits a  lower central series
\[ G \supseteq D:=(G,G) \supseteq (D,G) \supseteq \ldots  \]
The first quotient $ G/D $ is the homology group $ H_1(X,\Z)  $. We consider in this note 
 the second quotient $ D/(D,G)  $.  In particular when $ H_1(X,\Z)  $ is torsion free, we obtain a  description of $ D/(D,G)  $ in terms of the homology and cohomology of $ X $ (see Corollary \ref{tf} below).
 
 As an example, we recover in the last section the  isomorphism $D/(D,G)\cong \Z/2$ (due to Collino)
for the Fano surface  parametrizing the lines contained in a cubic threefold.

 \bigskip
 \section{The main result}
 
\begin{prop}
Let  $ X $  be a  connected space homotopic to a CW-complex, with $ H_1(X,\Z)  $ finitely generated.  Let $ G= \pi_1(X) $, $ D= \allowbreak (G,G)  $ its derived subgroup, $ \tilde{D} $ the subgroup of elements of $ G $ which are torsion in $ G/D $. 
The group $ D/(\tilde{D},G ) $ is canonically isomorphic to the cokernel of the map \[ \mu:  H_2(X,\Z) \rightarrow \mathrm{Alt}^2(H^1(X,\Z) ) \quad \mbox{given by }\   \mu(\sigma )(\alpha,\beta)= \sigma \frown (\alpha\wedge \beta) \ ,  \] where $ \mathrm{Alt}^2(H^1(X,\Z) ) $  is the group of skew-symmetric integral  bilinear forms on $ H^1(X,\Z)  $.
\end{prop}

\pr Let $ H $ be the quotient of $ H_1(X,\Z)  $ by its torsion subgroup; we put  $ V:=H\otimes_{\Z} \R $ and $ T:=V/H $. The quotient map $ \pi: V\rightarrow T $ is the universal covering of the real torus $ T $.

Consider the surjective homomorphism $\alpha: \pi_1(X)\rightarrow H  $. 
Since $ T $ is a $ K(H,1)  $, there is a continuous map $ a:X\rightarrow T $, well defined up to homotopy, inducing $ \alpha $ on the fundamental groups. Let $ \rho:X'\rightarrow X $ be the pull back by $ a $ of the \'etale covering $ \pi : V\rightarrow T$, so that $ X':=X\times_TV $ and $ \rho $ is the covering associated to the homomorphism $\alpha$.

Our key ingredient will be the map $ f: X\times V\rightarrow T $ defined by $ f(x,v)=a(x) -\pi(v)   $. It is a locally trivial fibration, with fibers isomorphic to $ X' $. Indeed the diagram 
\[ \xymatrix{X'\times V \ar[r]^{g}\ar[d]_{\mathrm{pr^{}_2}}& X\times V \ar[d]^f \\ V \ar[r]^\pi& T } \]
where $ g((x,v),w)=(x,v-w)   $, is cartesian.

It follows from this diagram that the monodromy action of $ \pi_1(T)=H  $ on $ H_1(X',\Z )  $ is induced by the action of $ H $ on $X' $; it is deduced from the action of $  \pi_1(X) $ on $ \pi_1(X' ) $ by conjugation in the exact sequence
\begin{equation}\label{sexpi}
 1\rightarrow \pi_1(X' )\qfl{\rho_*} \pi_1(X)\rightarrow H\rightarrow 1\ .  
 \end{equation}

The homology spectral sequence of the fibration $ f $ (see for instance \cite{Hu}) gives rise in low degree to a  five terms exact sequence
\begin{equation}\label{5}
H_2(X,\Z) \qfl{a_*} H_2(T,\Z) \longrightarrow H_1(X' ,\Z)_H \qfl{\rho_*}H_1(X,\Z) \longrightarrow H_1(T,\Z) \longrightarrow 0   \ ,
 \end{equation}
where $ H_1(X' ,\Z)_H $ denote the coinvariants of $ H_1(X' ,\Z) $ under the action of $ H $. 

The exact sequence (\ref{sexpi})  identifies  $ \pi_1(X') $ with $ \tilde{D} $, hence $ H_1(X' ,\Z) $ with 
$ \tilde{D}/(\tilde{D},\tilde{D})  $, the action of $ H $  being deduced from the action of $ G $ by conjugation. The group of coinvariants is the largest quotient of this group on which $ G $ acts trivially, that is, the quotient $  \tilde{D}/(\tilde{D},G) $. 

The exact sequence (\ref{5}) gives an isomorphism  $ \Ker \rho_*\iso\Coker a_* $. The map $ \rho_*: H_1(X' ,\Z)_H \rightarrow\allowbreak H_1(X,\Z) $ is identified with the natural map $ \tilde{D}/(\tilde{D},G) \rightarrow G/D $ deduced from the inclusions $\tilde{D}  \subset G $ and $ (\tilde{D},G) \subset D $. Therefore its kernel is $ D/(\tilde{D},G)  $. On the other hand
since $ T $ is a torus we have canonical isomorphisms 
\[ H_2(T,\Z) \iso \Hom (H^2(T,\Z),\Z)\iso   \mathrm{Alt}^2(H^1(T,\Z) ) \iso \mathrm{Alt}^2(H^1(X,\Z) ) \ , \]
through which $ a_* $ corresponds to $\mu  $, hence the Proposition.\qed

\begin{cor}\label{c1}
$1)$ There is a canonical surjective map $  D/(D,G) \rightarrow\Coker \mu $ with finite kernel.

$2)$ There are canonical exact sequences
\[ H_2(X,\Q) \qfl{\mu^{}_{\Q}}\mathrm{Alt}^2(H^1(X,\Q))\longrightarrow D/(D,G)\otimes \Q\rightarrow 0\] \[0\rightarrow \Hom(D/(D,G),\Q)\longrightarrow  \wedge^2H^1(X, \Q)\qfl{c^{}_\Q} H^2(X, \Q)\ ,\] where $c^{}_\Q$ is the cup-product map.
\end{cor}
\pr  
2) follows from 1), and from the fact that the transpose of $\mu^{}_\Q$ is $c^{}_\Q$. Therefore
in view of the Proposition, it suffices to prove that the kernel of the natural map  $ D/(D,G)\rightarrow  D/(\tilde{D} ,G)  $, that is, $ (\tilde{D} ,G)/(D ,G)  $, is finite. Consider the surjective homomorphism 
\[ G/D \otimes G/D \rightarrow D/(D,G) \]
deduced from $ (x,y) \mapsto xyx^{-1}y^{-1}  $. It maps $ \tilde{D} /D \otimes G/D $ onto $ (\tilde{D} ,G)/(D ,G) $; since $  \tilde{D} /D $ is finite and $ G/D $ finitely generated, the result follows.\qed

\begin{cor}\label{tf}
Assume that $ H_1(X,\Z)  $ is torsion free.

$ 1) $ The second quotient $ D/(D,G)  $ of the lower central series of $ G $  is canonically isomorphic to $ \Coker \mu $.

$ 2) $ For every ring  $ R $ the group 
$\Hom(D/(D,G), R)$  is canonically isomorphic to the kernel of the cup-product map  $c^{}_R: \wedge^2 H^1(X,R) \rightarrow H^2(X,R)$. 
\end{cor}
\pr 
We have $\tilde{D}=D  $ in that case, so $ 1) $ follows immediately from the Proposition. Since $ H_1(X,\Z)  $ is torsion free, the universal coefficient theorem provides an isomorphism $ H^2(X,R)\iso\allowbreak \Hom(H_2(X,\Z),R )   $, hence applying $ \Hom(-,R)  $ to the exact sequence
\[ H_2(X,\Z) \rightarrow  \mathrm{Alt}^2(H^1(X,\Z)) \rightarrow D/(D,G)\rightarrow 0   \]
gives $ 2) $.\qed

\begin{rem}
The Proposition and its Corollaries hold (with the same proofs) under weaker assumptions on $ X $, for instance for a connected space $ X $ which is paracompact, admits a universal cover and is such that $ H_1(X,\Z)  $ is finitely generated. We leave the details to the reader.
\end{rem}
\begin{rem}
For compact K\"ahler manifolds, the isomorphism $  \Hom(D/(D,G),\Q)\cong \Ker c^{}_{\Q} $ (Corollary \ref{c1}) is usually deduced from Sullivan's theory of minimal models  (see \cite{FG}, ch.\ 3); it can be used to prove that certain manifolds, for instance Lagrangian submanifolds of an abelian variety, have a non-abelian fundamental group.
\end{rem}
\section{Example: the Fano surface}

\medskip
Let $ V \subset \P^4$  be a smooth cubic threefold. The Fano surface $ F $ of $ V $ parametrizes the lines contained in $ V $. It is a smooth connected surface, which has been thoroughly studied in \cite{CG}. Its Albanese variety $A$ is canonically isomorphic to the intermediate Jacobian $ JV $ of $ V $, and the Albanese map $ a: F\rightarrow A $ is an embedding. Recall that $ A=JV $ carries a principal polarization
 $ \theta\in H^2(A,\Z)  $; for each integer $ k $ the class $ \dfrac{\theta^k}{k!} $ belongs to $ H^{2k}(A,\Z)  $. The class of $F$ in $H^6(A,\Z)$ is $\dfrac{\theta ^3}{3!}$ (\cite{CG}, Proposition 13.1).
 
 \begin{prop}
The maps $a^*:H^2(A,\Z)\rightarrow H^2(F,\Z)$ and $a_*: H_2(F,\Z)\rightarrow H_2(A,\Z)$ are injective and their images have index 2.
\end{prop}
 \pr We first recall that if $u:M\rightarrow N$ is a  homomorphism between two free $\Z$-modules of the same rank, the integer $|\det u|$ is well-defined : it is equal to the absolute value of the determinant of the matrix of $u$ for any choice of bases for $M$ and $N$. 
 If it is nonzero, it is equal to the index of $\im u$ in $N$.

 Poincar\'e duality identifies $a_*$ with the Gysin map $a_*:H^2(F,\Z)\rightarrow H^8(A,\Z)$, and also to the transpose of $a^*$. The composition $f : H^2(A,\Z)\qfl{a^*}H^2(F,\Z)\qfl{a_*}H^8(A,\Z)$ is the cup-product with the class $[F]=\dfrac{\theta ^3}{3!} $. We have $|\det a^*|=|\det a_*|\neq 0$  (\cite{CG}, 10.14), so it suffices to show that $|\det f|=4$.
 
 The principal polarization defines a unimodular skew-symmetric form on $H^1(A,\Z)$; we choose a symplectic basis $(\varepsilon _i,\delta _j)$ of  $H^1(A,\Z)$. Then 
 \[ \theta =\sum_i \varepsilon _i\wedge \delta _i \qquad \mbox{and} \qquad \dfrac{\theta ^3}{3!}=\sum_{i<j<k}(\varepsilon _i\wedge \delta _i)\wedge(\varepsilon _j\wedge \delta _j)\wedge(\varepsilon _k\wedge \delta _k) \ .\]
  If we identify by Poincar\'e duality $H^8(A,\Z)$ with the dual of $H^2(A,\Z)$, and $H^{10}(A,\Z)$ with $\Z$, $f$ is the homomorphism associated to the bilinear symmetric form $b:(\alpha ,\beta )\mapsto \alpha \wedge \beta \wedge \dfrac{\theta ^3}{3!}$, hence $|\det f|$ is the absolute value of the discriminant of $b$.
  Let us write
 $H^2(A,\Z)=M\oplus N$, where $M$ is spanned by the vectors $\varepsilon _i\wedge \varepsilon _j$, $\delta _i\wedge \delta _j$ and $\varepsilon _i\wedge \delta _j$ for $i\neq j$, and $N$ by the vectors $\varepsilon _i\wedge \delta _i$. 
 The decomposition is orthogonal with respect to $b$; 
 the restriction of $b$ to $M$ is unimodular, because the dual basis of $(\varepsilon _i\wedge \varepsilon _j,\delta _i\wedge \delta _j,\varepsilon _i\wedge \delta _j)$ is $(-\delta _i\wedge \delta _j,- \varepsilon _i\wedge \varepsilon _j,- \varepsilon _j\wedge \delta _i)$.  On $N$ the matrix of $b$ with respect to the basis $(\varepsilon _i\wedge \delta _i)$ is $E-I$, where $E$ is the 5-by-5 matrix with all  entries equal to 1. Since $E$ has rank $1$ we have $\wedge^kE=0$ for $k\geq 2$, hence
 \[\det(E-I)=-\det(I-E)=-I+\Tr E =4\ ;\]
 hence $|\det f|=4$.\qed
  
  \begin{cor}\label{fano}
Set $G=\pi _1(F)$ and $D=(G,G)$. The group $D/(D,G)$ is cyclic of order $2$.
\end{cor}
 Indeed $H_1(F,\Z)$ is torsion free \cite{C1}, hence the result follows from Corollary \ref{tf}.\qed
 
 \begin{rem}
The deeper topological study of  \cite{C1} gives actually the stronger result that $D$ is generated as a normal subgroup by an element $\sigma $ of order 2 (see \cite{C1}, and the correction in \cite{C2}, Remark 4.1). Since every conjugate of $\sigma $ is equivalent to $\sigma $ modulo $(D,G)$, this implies Corollary \ref{fano}.
\end{rem}
\begin{rem}
Choose a line $ \ell\in F $, and let $ C\subset F $ be the curve of lines incident to $ \ell $. Let $d: H^2(F,\Z)\rightarrow\allowbreak \Z/2$ be the homomorphism given by $ d(\alpha)= (\alpha\cdot [C])  $ $ (\mathrm{mod.}\ 2)  $. We claim that the image of  $a^*:H^2(A,\Z)\rightarrow H^2(F,\Z)$ is $\Ker d$. Indeed we have $ (C^2)=5 $ (the number of lines incident to two given skew lines on a cubic surface), hence $d([C])=1$, so that $\Ker d$ has index 2; thus it suffices to prove $d\rond a^*=0 $.
For $ \alpha\in H^2(A,\Z)  $, we have $d(a^*\alpha )=(a^*\alpha\cdot [C] ) = (\alpha\cdot a_*[C])\ \mathrm{mod.}\ 2$; this is 0 because the class $a_*[C]\in H^8(A,\Z)$ is equal to $2 \, \dfrac{\theta^4}{4!} $ (\cite{CG}, Lemma 11.5), hence  is divisible by 2.  

We can identify $a^*$  with the cup-product map $c$; thus we have an exact sequence
\begin{equation}
 0\rightarrow \wedge^2 H^1(F,\Z) \qfl{c} H^2(F,\Z) \qfl{d }
 \Z/2\rightarrow 0\qquad \mbox{with }\  d(\alpha)= (\alpha\cdot [C]) \  (\mathrm{mod.}\ 2)\, .  \end{equation}

\end{rem}

\end{document}